\documentclass[graybox]{svmult}

\usepackage{type1cm}        
\usepackage{makeidx}         
\usepackage[dvipdfmx]{graphicx}        
\usepackage{multicol}        
\usepackage[bottom]{footmisc}

\usepackage{newtxtext}       %
\usepackage[varvw]{newtxmath}       

\usepackage{natbib} 
\usepackage{booktabs}
\usepackage{lscape}
\usepackage{latexsym,indentfirst,ascmac}
\usepackage{hhline}
\usepackage{bm}

\makeindex       

\begin{document}

\title*{Semiparametric Estimation of Optimal Dividend Barrier for Spectrally Negative L\'{e}vy Process}
\titlerunning{Dividend Barrier Estimation for L\'{e}vy Process}
\author{Yasutaka Shimizu and Hiroshi Shiraishi}
\institute{Yasutaka Shimizu \at Waseda University, 3-4-1 Okubo, Shinjuku-ku, Tokyo, 169-8555, Japan, \email{shimizu@waseda.jp}
\and Hiroshi Shiraishi \at Keio University, 3-14-1 Hiyoshi, Kohoku-ku, Yokohama, Kanagawa, 223-8522, Japan, \email{shiraishi@math.keio.ac.jp}}
%
%
\maketitle

\abstract{We disucss a statistical estimation problem of an optimal dividend barrier when the surplus process follows a L\'{e}vy insurance risk process. 
The optimal dividend barrier is defined as the level of the barrier that maximizes the expectation of the present value of all dividend payments until ruin. 
In this paper, an estimatior of the expected present value of all dividend payments is defined based on ``quasi-process'' in which sample paths are generated by shuffling increments of a sample path of the L\'{e}vy insurance risk process. 
The consistency of the optimal dividend barrier estimator is shown. 
Moreover, our approach is examined numerically in the case of the compound Poisson risk model perturbed by diffusion.  
}

\section{Introduction}
In risk theory, surplus process is a very important model for understanding how the capital or surplus of an insurance company evolves over time. 
The classical model for the surplus process is the so-called ``Cram\'{e}r-Lundberg insurance risk model''. 
In this model, the insurance company collects premiums at a fixed rate $c>0$ from its customers. 
On the other hand, a customer can make a claim causing the surplus to jump downwards. 
The claim frequency follows a Poisson process, and the claim sizes are assumed to be  independent and identically distributed (i.i.d.).      
A natural generalization of the Cram\'{e}r-Lundberg model is a spectrally negative L\'{e}vy process also called ``L\'{e}vy insurance risk model'', which has been studied in many  actuarial literature, such as \cite{F11}, \cite{FS13}, \cite{K14} and so on. 
Thanks to the L\'{e}vy insurance risk model, we can grasp many realistic social phenomena in the surplus process such as the fluctuation of premium income, the effect of investment result, the effect of the small claim and so on.  
In this paper, we also suppose that the surplus of an insurance company follows the L\'{e}vy insurance risk process. 

In risk theory, the central topics are the ruin time or ruin probability, but there is a  dividend problem as one of the application.   
In the dividend problem introduced by \cite{D57} (especially, the so-called ``constant barrier strategy''), assuming that there is a horizontal barrier of level $\vartheta$, such that when an insurance company's surplus reaches level $\vartheta$, dividends are paid continuously such that the surplus stays at level $\vartheta$ until it becomes less than $\vartheta$. 
The optimal strategy is to maximize the expectation of the present value of all dividend payments and the ``optimal dividend barrier'' is defined as a barrier of level $\vartheta$ where the maximization can be achieved.     
\cite{B70}, \cite{L03}, \cite{G06}, \cite{L06}, \cite{L08} derived the optimal dividend barrier explicitly in some special models such as the Cram\'{e}r-Lundberg model with exponential claim amount distribution. 
On the other hand, \cite{K14} discussed a stochastic control problem for the optimal dividend strategy when the surplus process follows the L\'{e}vy insurance risk process.   
In these papers, the main concern is the property from a probabilistic point of view, but there is a limited contribution in the statistical point of view. 
From the statistical point of view in ruin theory, a ruin probability by \cite{C90} and \cite{S09},  a Gurber-Shiu function by \cite{FS13} and an optimal dividend problem in the Cram\'{e}r-Lundberg model by \cite{SL18} are discussed, respectively. 
In this paper, we discuss the statistical estimation problem in the L\'{e}vy insurance risk process. 

Considering the optimal dividend problem in a statistical estimation framework, it can be reduced to an M-estimation problem if we can define our optimal dividend barrier estimator as a maximizer of an objective function which corresponds to an estimator of the expectation of the present value of all dividend payments.  
Note that, for the usual M-estimator, the objective function sometimes called a contrast estimator is defined by a sample mean of i.i.d. random variables. 
In the same way, our contrast estimator would be defined by a sample mean of the present value of all dividend payments. 
However, since the present value of all dividend payments is path dependent, in order to construct the contrast estimator, we need to provide a number of independent copy of sample paths; consequently, it is impossible to observe multiple sample paths. 
In addition, in practical point of view, it is reasonable to assume that the surplus of a insurance portfolio is observable discretely not continuously, such as hourly, daily, monthly and so on. 
To overcome these problems, we introduce ``quasi-process'', that is, an approximation of the true L\'{e}vy insurance risk process. 
The quasi-process is composed by rearranging the increments of discretely observed data, thus, it is possible to generate multiple sample paths by changing the permutation. 
Essentially, it takes advantage of the exchangeability of the increments in the L\'{e}vy insurance risk process.   
Generating multiple quasi-process, it is possible to provide a number of (approximated) present value of all dividend payments, which implies that a contrast estimator can be defined.  
In our estimation procedure, the complexity of an estimator is characterized in the class of functions. 
In other words, our procedure is applicable not only to optimal dividend problem but also to many statistical inference problems defined as an M-estimation problem.    

This paper is organized as follows. 
Section 2 defines the surplus process following the L\'{e}vy insurance risk process and the true optimal dividend barrier as a maximizer of the expectation of the present value of all dividend payments. 
We define the quasi-process from discretely observed data and show its weak convergence in Section 3. 
Then, the optimal dividend barrier estimator is also defined.  
Section 4 shows the consistency of the optimal dividend barrier estimator. To do so, the uniform consistency for the contrast estimator in the function set is shown based on the empirical process theory. 
In Section 5, we examine our approach numerically. 
When the surplus process follows the compound Poisson risk model perturbed by diffusion discussed in \cite{L06}, it is numerically confirmed that our proposed estimator converges in probability to the true optimal dividend barrier as observe interval goes to $0$ and the size of permutation set goes to infinity. 
We place all the proofs of the theorems and lemmas in Section 6. 

\section{Optimal Dividend Barrier}
Given a stochastic basis $(\Omega,\mathcal{F},\mathbb{P};\mathbb{F})$ with a filtration $\mathbb{F}=(\mathcal{F}_t)_{t\ge 0}$, we consider a  $\mathbb{F}-$L\'{e}vy process $X=(X_t)_{t\ge 0}$ starting at $X_0=u$ of the form 
\begin{align}
X_t=u+ ct + \sigma W_t - S_t, \label{eq1}
\end{align}
where $u,\sigma \ge 0$, $c>0$, $W=(W_t)_{t\ge 0}$ is a Wiener process  and $S=(S_t)_{t\ge 0}$ is a pure-jump L\'{e}vy process, independent of $W$, with the characteristic exponent
\begin{align*}
\psi_S(\lambda) = 
\log \mathbb{E}[e^{i\lambda S_1}] = \int_{\mathbb{R}}\left(e^{i\lambda z} -1 - i\lambda z \bm{1}_{\{|z|\le 1\}} \right) \nu(dz). 
\end{align*}
When $\nu((-\infty,0))=0$ and $\int_{(0,\infty)}(1\wedge x)\nu(dx)<\infty$, $S$ is called a subordinator, that is, a special class of L\'{e}vy processes taking values in $[0,\infty)$ and having non-decreasing paths. 
Let $\mathbb{D}_{\infty}:=D[0,\infty)$ be a space of c\`{a}dl\`{a}g functions on $[0,\infty)$, 
and the subset $\tilde{\mathbb{D}}_{\infty}(\subset \mathbb{D}_{\infty})$ be also a space of c\`{a}dl\`{a}g functions on $[0,\infty)$, restricted as follows:

For all $X\in \tilde{\mathbb{D}}_{\infty}$
, $X$ has the form (\ref{eq1}), where $c>\mathbb{E}[S_1]$ and $S$ is a subordinator with 
$\int_{(0,\infty)}x^2\nu(dx)<\infty$.     

\noindent Then, the $\tilde{\mathbb{D}}_{\infty}$ belongs to the class of spectrally negative L\'{e}vy processes with $\mathbb{E}[S_t^2]<\infty$ for all $t>0$ (see, e.g., \cite{S99}).  
In this paper, we suppose that $X\in \tilde{\mathbb{D}}_{\infty}$ is an insurance risk process where $u(=X_0)$ is the insurer's initial surplus, $c$ is a given premium rate per unit time, with the net profit condition $c>\mathbb{E}[S_1]$, and $S=(S_t)_{t\ge 0}$ is the aggregate claims process. 

Let $\Theta=(u,\bar{\vartheta})\subset \mathbb{R}$, where $\bar{\vartheta} 
$ is a known positive value. 
For an insurance risk process $X\in\tilde{\mathbb{D}}_{\infty}$ and   $\vartheta\in \Theta$, we introduce a process $\xi^{\vartheta}=(\xi_t^{\vartheta})_{t\ge 0}$ by $\xi^{\vartheta}_0=0$ and 
$$\xi^{\vartheta}_t = \vartheta \vee \bar{X}_t - \vartheta = \left(\bar{X}_t - \vartheta \right)\vee 0\quad \mathrm{for}\ t>0,$$ 
where $\bar{X}_t=\sup_{0\le s< t}X_s$. 
Note that $\xi^{\vartheta}\equiv \xi^{\vartheta}(X)$ is called the  dividend strategy consisting of a process with initial value zero, which has paths that are left-continuous, non-negative, non-decreasing and adapted to the filtration of insurance risk process $X$ defined by (\ref{eq1}). 
Let $\Xi \equiv \Xi(X)=\left\{\xi^{\vartheta}(X)|\vartheta\in\Theta \right\}$ be the family of dividend strategies, and for each $\xi^{\vartheta}\in\Xi$, write $\tau^{\vartheta}\equiv \tau^{\vartheta}(X) = \inf\{t>0|U_t^{\vartheta}:=X_t-\xi_t^{\vartheta}<0\}$ for the time of ruin under the dividend strategy $\xi^{\vartheta}$. 
Here we call $U^{\vartheta}=(U_t^{\vartheta})_{t\ge 0}$ the controlled risk process and $\tau^{\vartheta}$ the time of ruin for the controlled risk process (see, e.g., 
\cite{K14}); $\xi_t^{\vartheta}$ represents the cumulative dividends that the insurer has paid out until the time $t$ under a dividend strategy which the dividend payments are continued while the controlled risk process attains $\vartheta$ up to the time of ruin $\tau^{\vartheta}$ (see, e.g., \cite{L08}). 
The expected present value of all dividend payments, with discounting at rate $r> 0$, associated with the dividend strategy $\xi^{\vartheta}$ is given by 
\begin{align}
v(\xi^{\vartheta}) = \mathbb{E}\left[ h^{\vartheta}(X) \right] ,\quad h^{\vartheta} (X) = \int^{\tau^{\vartheta}(X)}_{0} e^{-rt}d\xi_t^{\vartheta}(X).\label{eq3}
\end{align}
\cite{L08} and \cite{Y15} discussed the concavity for $v(\xi^{\vartheta})$ 
under some conditions. 
We suppose that $v(\xi^{\vartheta})$ is a  bounded, infinitely differentiable, and strictrly concave function with respect to $\vartheta\in\Theta$. 
Then, for any $\epsilon>0$, there exists $\vartheta_0\in \Theta$ such that for all $\vartheta\in \Theta$ satisfying $|\vartheta -\vartheta_0|>\epsilon$, it follows $v(\xi^{\vartheta})<v(\xi^{\vartheta}_0)$. 
We assume the proper proerty for $v(\xi^{\vartheta})$ in our main theorem (Theorem \ref{th2}).  
\color{black} 
In insurance risk theory, the expected present value of a \textit{ruin-related `loss' up to time of ruin} 
is often discussed (e.g., \cite{F11} and \cite{FS13}). 
Among them, 
the dividend problem discussed in \cite{D57} consists of solving the stochastic control problem $v(\xi_*):= \sup_{\xi^{\vartheta} \in \Xi} \mathbb{E}\left[ h^{\vartheta}(X) \right]$ which corresponds to a optimization problem
\begin{align}
\vartheta_0 := \arg\max_{\vartheta \in\bar{\Theta}}\mathbb{E}\left[ h^{\vartheta}(X) \right].\label{eq4}
\end{align}
In this paper, we consider statistical estimation problem for $\vartheta_0$ when we observe an insurance risk process $X\in \tilde{\mathbb{D}}_{\infty}$ discretely. 

\section{Estimation of Optimal Dividend Barrier}
Let $\mathcal{D}_{\infty}$ be the Borel field on $\mathbb{D}_{\infty}$ generated by the Skorokhod topology. 
We denote a distribution of $X$ on $\mathcal{D}_{\infty}$ by $P:=\mathbb{P}\circ X^{-1}$ and write 
\begin{align*}
Pf := \int_{\mathbb{D}_{\infty}} f(x) P(dx) = \mathbb{E}[f(X)],
\end{align*}  
for a measurable function $f:\mathbb{D}_{\infty} \to \mathbb{R}$. 
Suppose that for a $B \in \mathbb N$, random elements $X^{(1)},X^{(2)},\ldots,X^{(B)}$ are independent copies of process $X\in 
\tilde{\mathbb{D}}_{\infty}(\subset \mathbb{D}_{\infty})$, 
and denote its empricial measure 
as 
\begin{align*}
\mathbb{P}_{B}^* := \frac{1}{B}\sum_{\beta =1}^{B} \delta_{X^{(\beta)}},
\end{align*}  
where $\delta_x$ is the delta measure concentrated on $x\in \tilde{\mathbb{D}}_{\infty}$. 
In practice, it is often impossible to observe the independent copies of $X$ and to observe the sample path continuously. 
To overcome these problems, we consider a construction of ``multiple quasi-processes'' from a discrete sample path. 
Suppose that we observe a discrete sample path from a insurance risk process $X=(X_t)_{t\ge 0}\in \tilde{\mathbb{D}}_{\infty}$, where the discrete sample path consists of $\{X_{t_k}\}_{k=0,1,\ldots,n}$ with 
\begin{align*}
0=t_0<t_1<\cdots <t_n=T,\quad h_n\equiv t_k-t_{k-1}.
\end{align*}  
Let $\mathbb{X}=(\Delta_1X,\Delta_2X,\ldots,\Delta_nX)$ be a vector of increments with $\Delta_kX:=X_{t_k}-X_{t_{k-1}}$, and let
\begin{align*}
\Lambda_n:=\left\{i_m=\left(
  \begin{array}{cccc}
   1&2&\cdots&n    \\
   i_m(1)&i_m(2)&\cdots &i_m(n)    \\
  \end{array}
\right)\Bigr| m=1,2,\ldots,n!
  \right\}
\end{align*}
be a family of all the permutations of $(1,2,\ldots,n)$. 
Since $\Delta_kX$, $1\le k\le n$, are i.i.d. for each $n$, $\mathbb{X}$ is \textit{exchangeable}, that is, for any permutation $i\in\Lambda_n$, 
\begin{align*}
i(\mathbb{X}):=(\Delta_{i(1)}X,\ldots,\Delta_{i(n)}X)
\end{align*}  
has the same distribution as $\mathbb{X}$.
\begin{definition}\label{def1}
For given $\mathbb{X}$ and $i\in\Lambda_n$, a stochastic process $\hat{X}^{i,n}=(\hat{X}_t^{i,n})_{t\ge 0}$ given by
\begin{align*}
\hat{X}_t^{i,n} = u + \sum_{k=1}^n \Delta_{i(k)}X \cdot \bm{1}_{[t_k,\infty)}(t)
\end{align*}
is said to be a \textit{quasi-process} of $X$ for a permutation $i\in\Lambda_n$.   
\end{definition}

Note that a path of the quasi-process $\hat{X}^{i,n}=(\hat{X}^{i,n}_t)_{t\ge 0}$ belongs to $\mathbb{D}_{\infty}$ (but not to $\tilde{\mathbb{D}}_{\infty}$), a right continuous step function that has a jump at $t=t_k$ ($k=1,2,\ldots,n$) with the amplitude $\Delta_{i(k)}X$. 
For the discrete sampling scheme, we impose the following assumption.\\
\ \\
\textbf{\textit{Assumption}} $\bm{1}$ (High-Frequency sampling in the Long Term; HFLT)
\begin{flushleft}
\quad $h_n\to 0$ and $T=nh_n\to \infty$ as $n\to\infty$. 
\end{flushleft}

\cite{S22} showed the followings under HFLT.

\newpage
\begin{theorem}
Under Assumption 1, we have, for any sequence of permutations $\{i^n\}\subset \Lambda_n$,
\begin{align*}
\hat{X}^{i^n,n}\leadsto X\quad \mathrm{in}\ \mathbb{D}_{\infty}\quad \mathrm{as}\ n\to \infty. 
\end{align*}
\end{theorem}
For a given size $\alpha_n(\le n!)$, let $A_n:=\{i_{(1)},\ldots, i_{(\alpha_n)}\}$ be a set of i.i.d. samples drawn uniformly from $\Lambda_n$, i.e., for a given $m=1,2,\ldots, n!$,
$$ \mathbb P(i_{(k)}=i_m) =\frac{1}{n!}\quad \mathrm{for\ every\ } k=1,2,\ldots, n!.$$
Based on $A_n$, we introduce two empirical measures $\mathbb{P}^*_{\alpha_n}$ and $\mathbb{P}_{\alpha_n}$ by   
\begin{align*}
\mathbb{P}^*_{\alpha_n} := \frac{1}{\alpha_n}\sum_{k=1}^{\alpha_n} \delta_{X^{(k)}},\quad 
\mathbb{P}_{\alpha_n} := \frac{1}{\alpha_n}\sum_{k=1}^{\alpha_n} \delta_{\hat{X}^{i_{(k)},n}}.
\end{align*}
Next, we propose an estimator of $\vartheta_0$ defined by (\ref{eq4}) based on the empirical measure of the quasi-process. 
\begin{definition}
Given a vector of increments $\mathbb{X}$ and permutation sets $\{A_n\}_{n\in\mathbb{N}}$, we denote a maximum contrast estimator of $\vartheta_0$ defined by (\ref{eq4}) as
\begin{align*}
\hat{\vartheta}_n = \mathrm{arg}\max_{\vartheta\in\bar{\Theta}} \mathbb{P}_{\alpha_n}h^{\vartheta}, 
\end{align*}
where
\begin{align*}
\mathbb{P}_{\alpha_n}h^{\vartheta} = \frac{1}{\alpha_n}\sum_{i\in A_n} h^{\vartheta}(\hat{X}^{i,n})= \frac{1}{\alpha_n}\sum_{i\in A_n}\int^{\tau^{\vartheta}(\hat{X}^{i,n})}_{0}e^{-rt}d\xi_t^{\vartheta}(\hat{X}^{i,n}).
\end{align*}
\end{definition}
For the moments of $h^{\vartheta}(\hat{X}^{i,n})$, we have following result. 
\begin{lemma}
Under Assumption 1, we have, for any $\vartheta\in\Theta$ and $i\in \Lambda_n$, 
\begin{align*}
\mathbb{E}\left[h^{\vartheta}(\hat{X}^{i,n})^m\right]=O(1),\quad m=1,2.
\end{align*}
\end{lemma}
\section{Asymptotic Results}
Our main result in this paper is to provide the consistency for $\hat{\vartheta}_n$ defined in Definition 2. 
To do so, we assume that a size of permutation sets $\alpha_n:=\sharp A_n$ satisfies followings. \\
\ \\
\textbf{\textit{Assumption}} $\bm{2}$ (Size of permutation sets)
\begin{flushleft}
\quad $\displaystyle \frac{n}{\alpha_n}\to 0$ as $n\to\infty$. 
\end{flushleft}

We recall that the empirical measure of the quasi-process $\mathbb{P}_{\alpha_n}$ is asymptotically equivalent in law with $\mathbb{P}^*_{\alpha_n}$ based on the independent copy $X^{(1)},\ldots,X^{(\alpha_n)}$ of the process $X$. 
Moreover, we introduce a sequence of a family of measurable functions $\mathcal{H}=\{\mathcal{H}_n\}_{n\in \mathbb{N}}$ on $\tilde{\mathbb{D}}_{\infty}$, where $\mathcal{H}_n$ is a family of measurable functions $h_n^{\vartheta}:\tilde{\mathbb{D}}_{\infty} \to \mathbb{R}$ for each $\vartheta \in \Theta$, given by
\begin{align}
h_n^{\vartheta}(X) = \int^{\tau_n^{\vartheta}(X)}_{0} e^{-rt}d\xi_{n,t}^{\vartheta}(X), \label{eq7}
\end{align}
where $\tau_n^{\vartheta}(X)=\tau^{\vartheta}(\hat{X}^{i_{\mathrm{id}},n})$ and $\xi_{n,t}^{\vartheta}(X)=\xi_{t}^{\vartheta}(\hat{X}^{i_{\mathrm{id}},n})$ for all $X\in \tilde{\mathbb{D}}_{\infty},\vartheta\in\Theta$ and $n\in\mathbb{N}$. 
Here $i_{\mathrm{id}}\in \Lambda_n$ is an identical permutation, i.e., $i_{\mathrm{id}}=\left(
  \begin{array}{cccc}
   1&2&\cdots&n    \\
   1&2&\cdots &n   \\
  \end{array}
\right)$.  
For the class $\mathcal{H}_n=\{h_n^{\vartheta}:\tilde{\mathbb{D}}_{\infty}\to \mathbb{R} |\vartheta\in\Theta\}$, we denote by $N(\epsilon,\mathcal{H}_n,L^1(\mathbb{P}_{\alpha_n}))$ the covering number of $L^1(\mathbb{P}_{\alpha_n})$ which is the minimum number of $\epsilon$-balls needed to cover $\mathcal{H}_n$, where an $\epsilon$-ball around a function $g\in L^1(\mathbb{P}_{\alpha_n})$ being the set $\{h_n^{\vartheta}\in L^1(\mathbb{P}_{\alpha_n})|\ \|h_n^{\vartheta}-g\|_{\mathbb{P}_{\alpha_n},1}=\mathbb{P}_{\alpha_n}(|h_n^{\vartheta}-g|) <\epsilon\}$, with $\|\cdot \|_{\mathbb{P}_{\alpha_n},1}$ being the $L^1(\mathbb{P}_{\alpha_n})$-norm. 
In addition, we denote by $N_{[]}(\epsilon,\mathcal{H}_n,L^1(\mathbb{P}_{\alpha_n}))$ the bracketing number which is the minimum number of $\epsilon$-brackets in $L^1(\mathbb{P}_{\alpha_n})$ needed to ensure that every $h_n^{\vartheta}\in \mathcal{H}_n$ lines in at least one bracket, where an $\epsilon$-bracket in $L^1(\mathbb{P}_{\alpha_n})$ is a pair of functions $l,u\in L^1(\mathbb{P}_{\alpha_n})$ with $\mathbb{P}_{\alpha_n}\left(\bm{1}_{\{l(X)\le u(X)\}}\right)=1$ and $\|l-u\|_{\mathbb{P}_{\alpha_n},1}\le \epsilon$. 
For the covering number and bracketing number, we have following result. 
\begin{lemma}
Under Assumptions 1 and 2, we have, for any $\epsilon>0$
\begin{align*}
\mathbb{E}\left[N(\epsilon,\mathcal{H}_n,L^1(\mathbb{P}_{\alpha_n}))\right]= o(\alpha_n)\quad \mathrm{and}\quad \mathbb{E}\left[N_{[]}(\epsilon,\mathcal{H}_n,L^1(\mathbb{P}_{\alpha_n}))\right]= o(\alpha_n).
\end{align*}
\end{lemma}

This lemma implies that both of the covering number and bracketing number diverge slower than $\alpha_n$. 
This result is applied for the proof of uniformly consistency below. 
The following result is due to a slight modification by Kosorok (2008; Theorem 8.15).
\begin{lemma}
Under Assumptions 1 and 2, we have
\begin{align*}
\sup_{\vartheta\in \bar{\Theta}}|(\mathbb{P}_{\alpha_n}-P)h_n^{\vartheta}|=: \|\mathbb{P}_{\alpha_n}-P\|_{\mathcal{H}_n}\stackrel{p}{\to}0.
\end{align*}
\end{lemma}

This lemma shows that two measure $\mathbb{P}_{\alpha_n}$ and $P$ are asymptotically equivalent on  the function space $\mathcal{H}_n$. 
On the other hand, the true optimal dividend barrier $\vartheta_0$ defined by (\ref{eq4}) is a maximizer of $\mathbb{E}\left[h^{\vartheta}(X)\right]= Ph^{\vartheta}$, where $h^{\vartheta}\notin \mathcal{H}_n$. 
Hence, we have to evaluate the difference between $h^{\vartheta}$ and $h_n^{\vartheta}\in \mathcal{H}_n$ defined by (\ref{eq3}) and (\ref{eq7}) based on the measure $P$. 
The following lemma is also applied for the proof of our main result.  
\begin{lemma}
Under Assumptions 1 and 2, we have
\begin{align*}
\sup_{\vartheta \in \bar{\Theta}}|P(h_n^{\vartheta}-h^{\vartheta})|\to 0.
\end{align*}
\end{lemma}

By Lemmas 3 and 4, we can show the consistency for $\hat{\vartheta}_n$, as follows.  
\begin{theorem}\label{th2}
Suppose that Assumptions 1 and 2 are hold, and that there exists $\vartheta_0\in\Theta$ such that, for any $\epsilon>0$,
\begin{align}
\sup_{\vartheta\in\bar{\Theta}:|\vartheta -\vartheta_0|>\epsilon} Ph^{\vartheta}< Ph^{\vartheta_0}.\label{eq8}
\end{align}
\color{black}
Then, $\hat{\vartheta}_n$ is weakly consistent to $\vartheta_0$, i.e., 
\begin{align*}
\hat{\vartheta}_n \stackrel{p}{\to} \vartheta_0,\quad n\to\infty.
\end{align*}
\end{theorem}
%
\section{Numerical Results}
In this section, we present simulation results to evaluate the finite-sample performance of the proposed estimator of the optimal dividend barrier based on the discrete sample from spectrally negative L\'{e}vy processes. 
We consider the following data generating process (DGP), sampling scheme, permutation set and discount rate. 
\begin{itemize}
  \item \textbf{DGP} (\textbf{Brownian motion} $\bm{+}$ \textbf{compound Poisson process}): Let $X_t=u+ ct+\sigma W_t - S_t$, where $u=10,\ c=15,\ \sigma=2$, $W=(W_t)_{t\ge 0}$ is a standard Brownian motion, and $S_t=\sum_{r=1}^{N_t}\xi_r$ is a compound Poisson process. The Poisson process $N=(N_t)_{t\ge 0}$ has intensity $\lambda>0$. We set $\lambda=5$. The jump size $\{\xi_r\}$ is a sequence of i.i.d. random variables having exponential distribution with parameter 
$1/2$, that is, $\mathbb{E}[\xi_r]=2$.\color{black}
  \item \textbf{Sampling scheme}: We consider the sampling interval $h_n=1,0.1,0.01,0.001$ and the terminal $T=100$ (fixed), which implies that sample size $n$ is $n=T/h_n=100,1000,10000$, respectively. 
  \item \textbf{Permutation set}: We consider the subset of permutation set $A_n=\{i_{m_j}|j=1,\ldots,\alpha_n\}\subset \Lambda_n$ with $\alpha_n=10,100,1000$, where the suffix $m_j$ is independently selected with same probability from $\{1,2,\ldots,n!\}$. 
  \item \textbf{Discount rate}: We set $r=0.2$. 
\end{itemize}
\subsection{Quasi Process} 
\begin{figure}[htbp]
  \begin{center}
    \includegraphics[width=100mm,height=50mm]{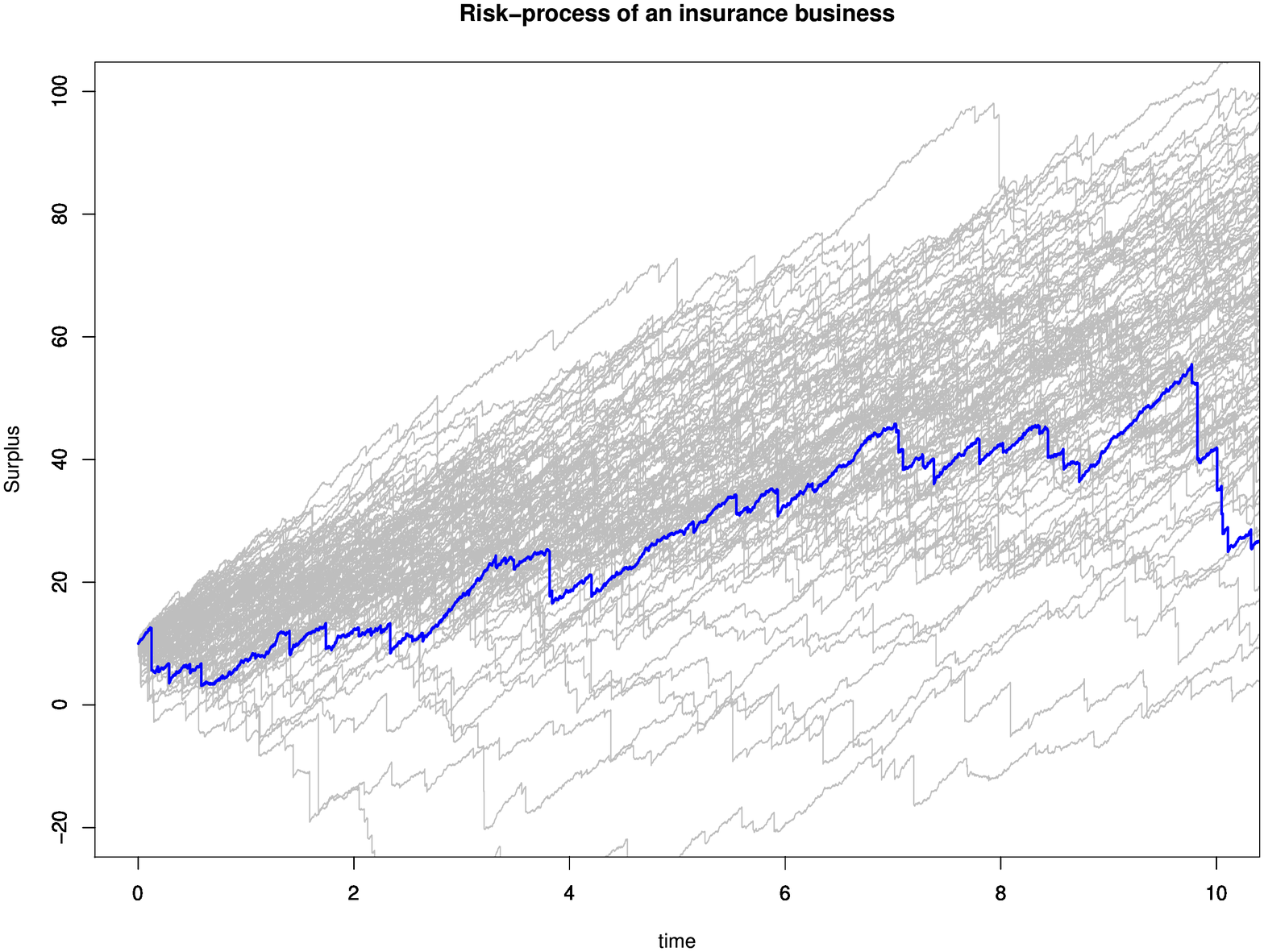}
  \end{center}
  \caption{100 sample paths for the risk process of an insurance business $X$.}
  \label{fig1}
\end{figure}

\begin{figure}[htbp]
  \begin{center}
    \includegraphics[width=70mm,height=50mm]{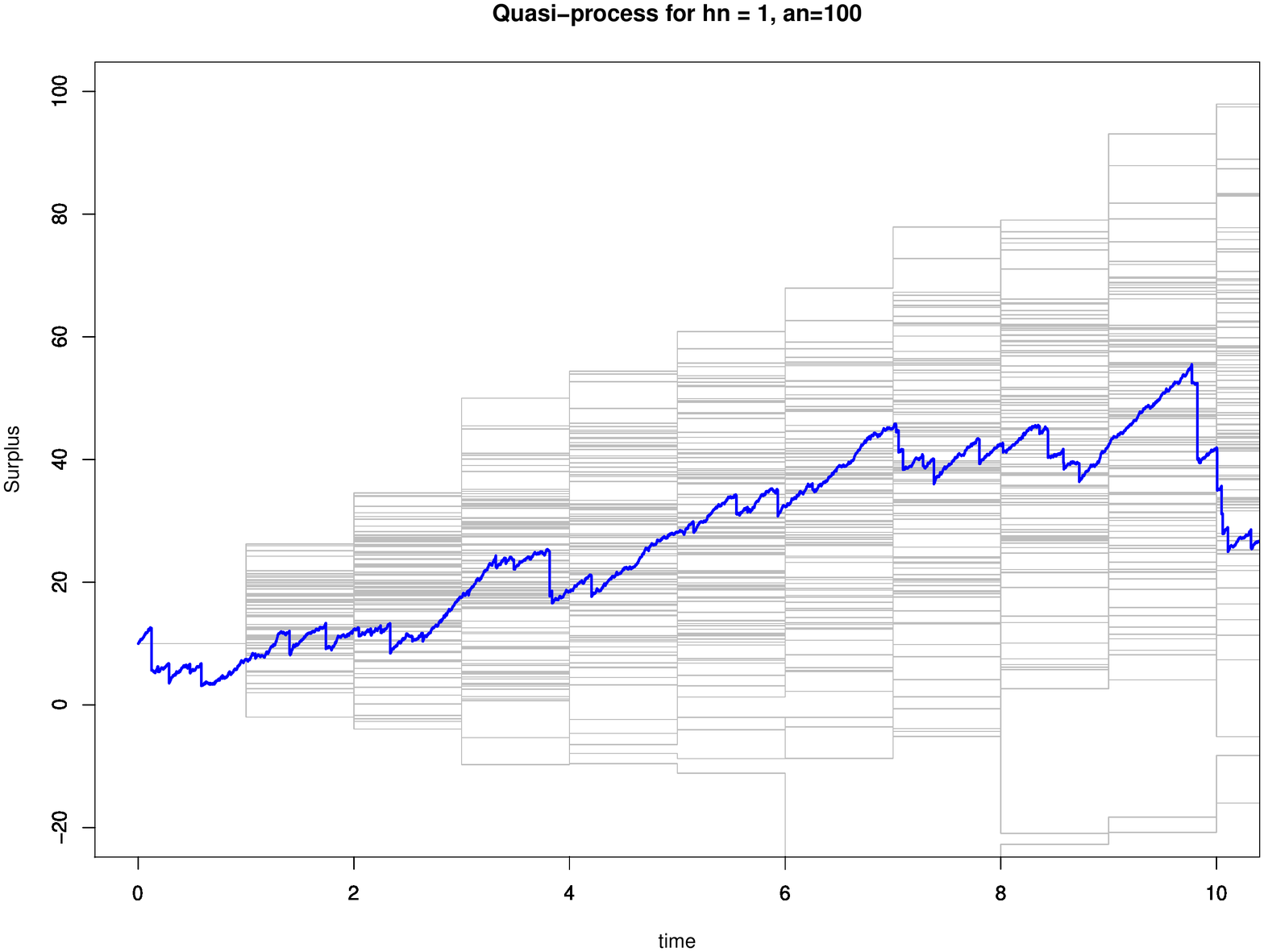}
    \includegraphics[width=70mm,height=50mm]{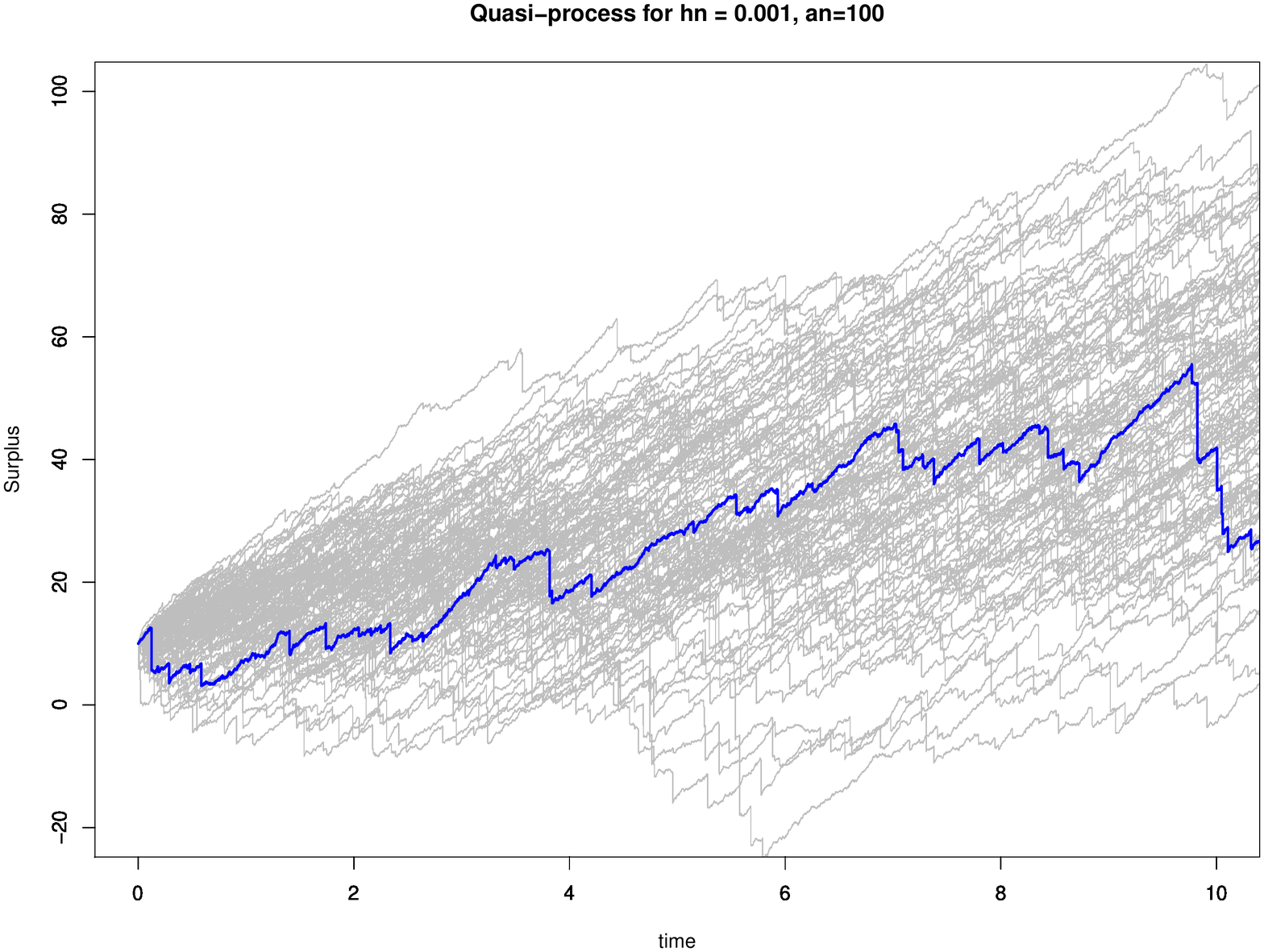}
  \end{center}
  \caption{(Discretely) observed sample path (blue line) and $\alpha_n(=100)$ sample paths for the quasi-process (top ($h_n=1$) and bottom ($h_n=0.001$)).}
  \label{fig2}
\end{figure}

We first examine the finite sample performance of the quasi-process $\hat{X}^{i,n}=(\hat{X}_t^{i,n})_{t\ge 0}$ for each $h_n$ and $\alpha_n$. 
Figure \ref{fig1} shows 100 sample paths for the risk process of an insurance business $X=(X_t)_{t\ge 0}$ defined above. 
It looks that we can not know the distribution of $X$ only from one sample path without any additional assumption. 
In this study, we consider such a situation. 
When we suppose that only one sample path is observed discretely, we would like to know its  distribution. 
In Figure \ref{fig1}, the blue line is observed discretely. 
Under the sampling scheme defined above, we can construct a number of sample paths of the quasi-process $\hat{X}^{i,n}$ from one sample path. 
Then, we can approximate the distribution of $X$ based on these sample paths. 
In Figure \ref{fig2}, the blue line is a observed (but discretely) sample path from the stochastic process $X$ (this is the same as Figure \ref{fig1}). 
From this sample path, we construct $\alpha_n$ sample paths of the quasi-process $\hat{X}^{i,n}$ based on Definition \ref{def1}. 
Each sample path depends on the observed sample path and the permutation $i\in A_n \subset \Lambda_n$. 
The top figure shows the case of $h_n=1$ and $\alpha_n=100$, and the bottom figure shows the case of  $h_n=0.001$ and $\alpha_n=100$. 
It looks that the top figure is not, but the bottom figure is well approximated the distribution of $X$. 
This phenomenon comes from the exchangeability of the increments of the L\'{e}vy processes and if the sampling interval $h_n$ is sufficiently small and the size of permutation set $\alpha_n$ is sufficiently large, we can well approximate the distribution of $X$ even from only one sample path.   

\subsection{Maximum Contrast Estimator} 
\begin{figure}[htbp]
  \begin{center}
    \includegraphics[width=100mm,height=50mm]{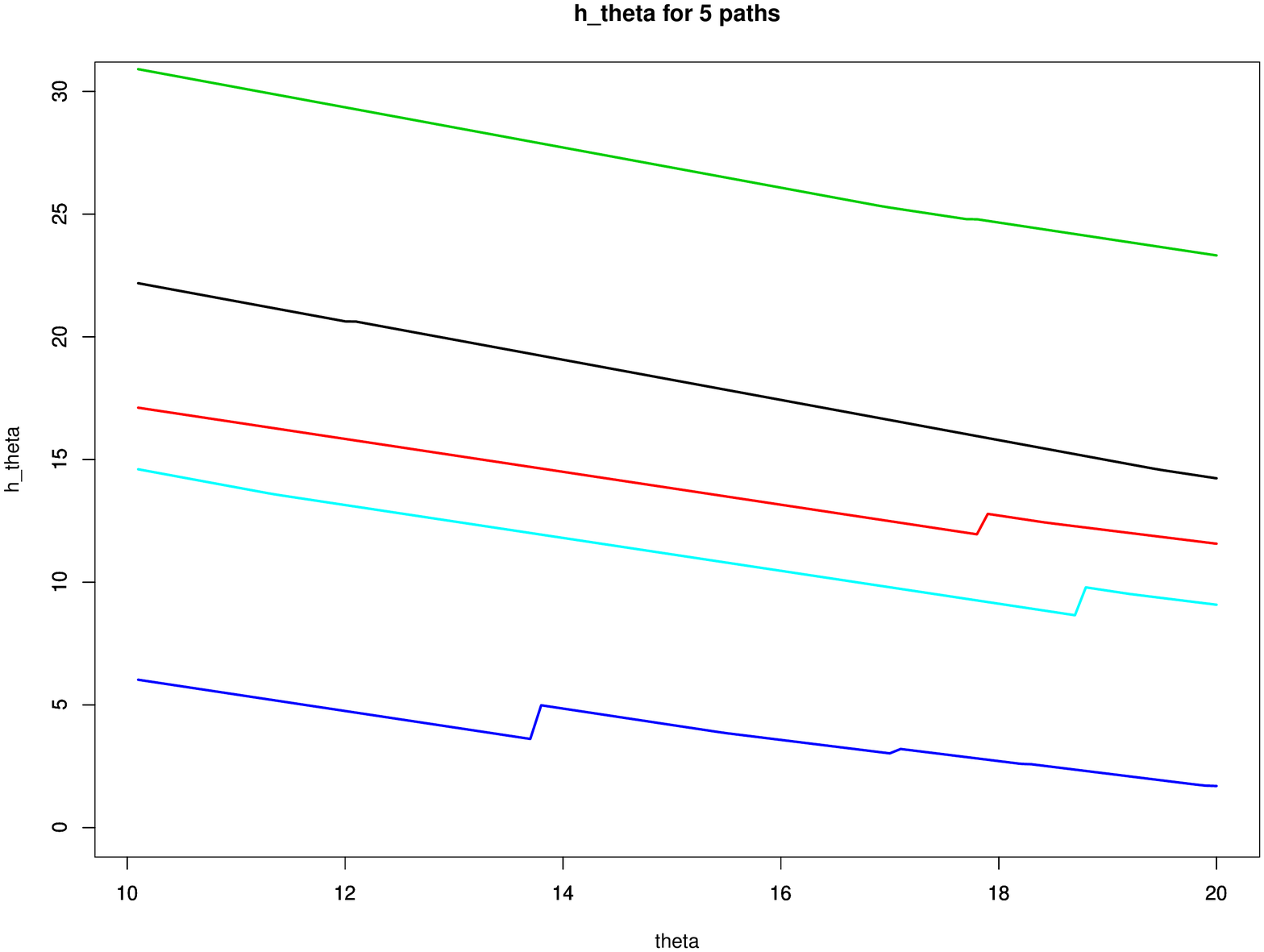}
  \end{center}
  \caption{Plots of $h_n^{\vartheta}$ for 5 sample paths of the quasi-process.}
  \label{fig3}
\end{figure}

Next, we examine the behavior of the objective function $h_n^{\vartheta}(\hat{X}^{i,n})$. 
Given a sample path of the quasi-process $\hat{X}^{i,n}$, we can construct $h_n^{\vartheta}(\hat{X}^{i,n})$ as a function of the parameter $\vartheta$. 
Since each sample path of the quasi-process $\hat{X}^{i,n}=(\hat{X}^{i,n})_{t\ge 0}$ is locally constant on time $t$, we can write
\begin{align*}
h_n^{\vartheta}(\hat{X}^{i,n})= \int^{\tau_n^{\vartheta}(\hat{X}^{i,n})}_{0} e^{-rt} d\xi_{n,t}^{\vartheta}(\hat{X}^{i,n}) =  \sum_{k=1}^n \bm{1}_{\{\tau_n^{\vartheta}(\hat{X}^{i,n})>t_k\}} e^{-rt_k} \Delta_{k}\xi_{n}^{\vartheta}(\hat{X}^{i,n}),
\end{align*}
where $\Delta_{k}\xi_{n}^{\vartheta}(\hat{X}^{i,n}) = \xi_{n,t_k}^{\vartheta}(\hat{X}^{i,n}) -\xi_{n,t_{k-1}}^{\vartheta}(\hat{X}^{i,n})$ with $\xi_{n,t_0}^{\vartheta}(\hat{X}^{i,n}) =0$.  
Figure \ref{fig3} shows the plots of $h_n^{\vartheta}(\hat{X}^{i,n})$ for five sample paths of the quasi-process $\hat{X}^{i,n}$. 
In this figure, the horizontal axis represents the magnitude of $\vartheta$ and the vertical axis represents the magnitude of $h_n^{\vartheta}$.   
Under a fixed sample path, it can be seen that $h_n^{\vartheta}$ is a locally decreasing function with some positive jumps. 
The locally decreasing property is that the total dividend amount tends to decrease as the dividend barrier $\vartheta$ increases while the ruin time $\tau_n^{\vartheta}$ is fixed. 
On the other hand, the existence of positive jump shows that the total dividend amount increases discontinuously since the ruin time is extended at some $\vartheta$.    

\begin{figure}[htbp]
  \begin{center}
    \includegraphics[width=110mm,height=50mm]{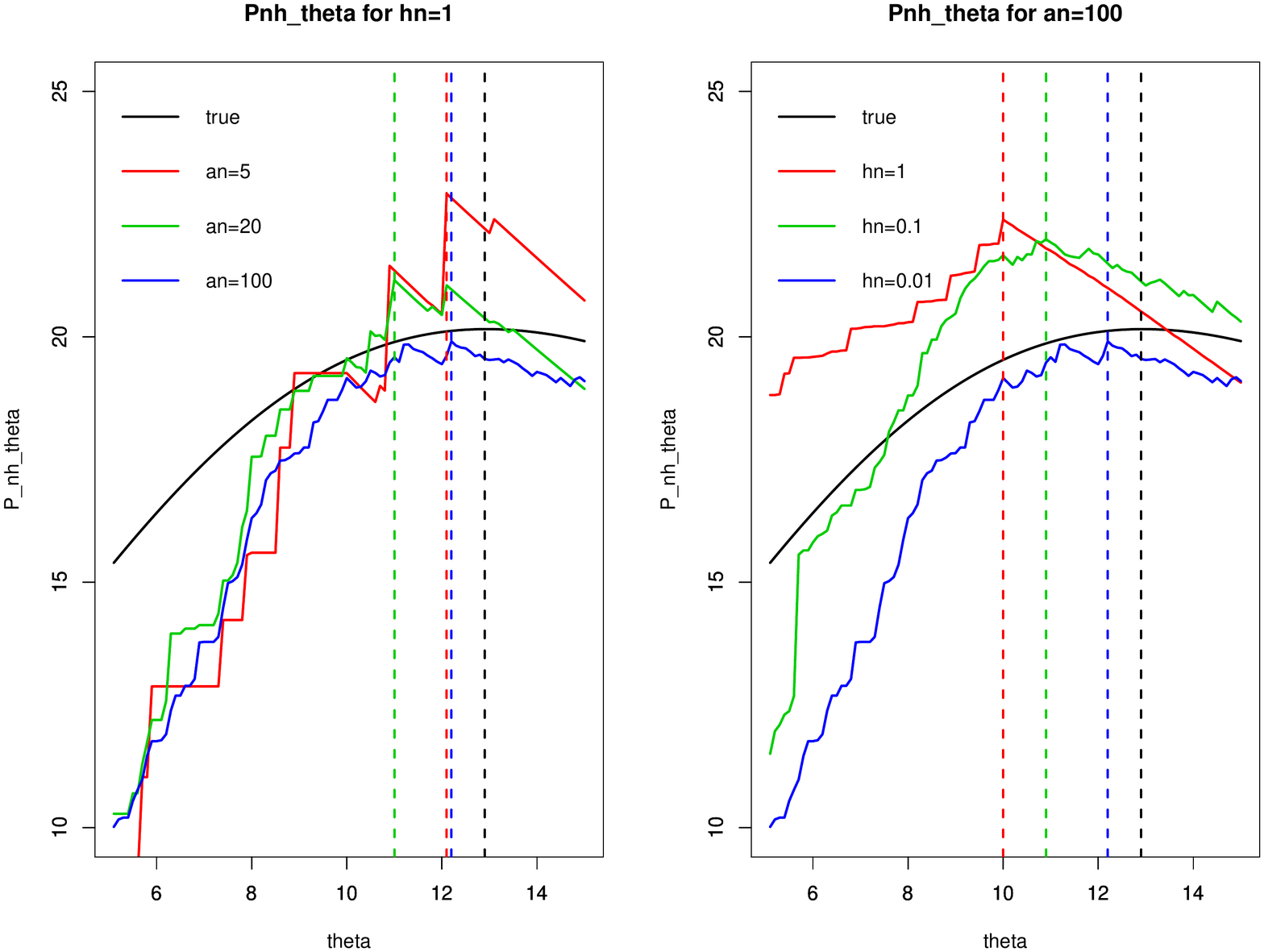}
  \end{center}
  \caption{Plots of the contrast function for the size of permutation set $\alpha_n=5,20,100$ under fixed sampling interval $h_n=1$ (left), and for $h_n=1,0.1,0.01$ under fixed $\alpha_n=100$ (right). The dotted line shows these maximization points, and the black line represents the true objective function $\mathbb{E}[h^{\vartheta}(X)]$.}
  \label{fig4}
\end{figure}

Figure \ref{fig4} shows the behavior of the contrast function $\mathbb{P}_{\alpha_n}h_n^{\vartheta}= \frac{1}{\alpha_n}\sum_{i\in A_n}h_n^{\vartheta}(\hat{X}^{i,n})$
 for some $h_n$ and $\alpha_n$. 
The left figure shows plots of the contrast function for the size of permutation set $\alpha_n=5,20,100$ under fixed sampling interval $h_n=1$, and the dotted line shows its maximization point. 
It can be seen that the function approaches the true function as $\alpha_n$ increases which implies that the maxmization point tends to the true maximation point, that is, our proposed estimator $\hat{\vartheta}_n$ converges to the true optimal dividend barrier $\vartheta_0$.  
On the other hand, the right figure shows plots of the contrast function for $h_n=1,0.1,0.01$ under fixed $\alpha_n=100$. 
It can be seen that the function approaches the true function as $h_n$ decreases which implies that our estimator converges to the true optimal dividend barrier.  
In both figures, the black line represents the true objective function $\mathbb{E}[h^{\vartheta}(X)]$ (see, e.g., \cite{L06}). 
These figures confirm the validity of the theoretical result in Theorem \ref{th2}. 

\subsection{Simulation result} 
\begin{figure}[htbp]
  \begin{center}
    \includegraphics[width=110mm,height=50mm]{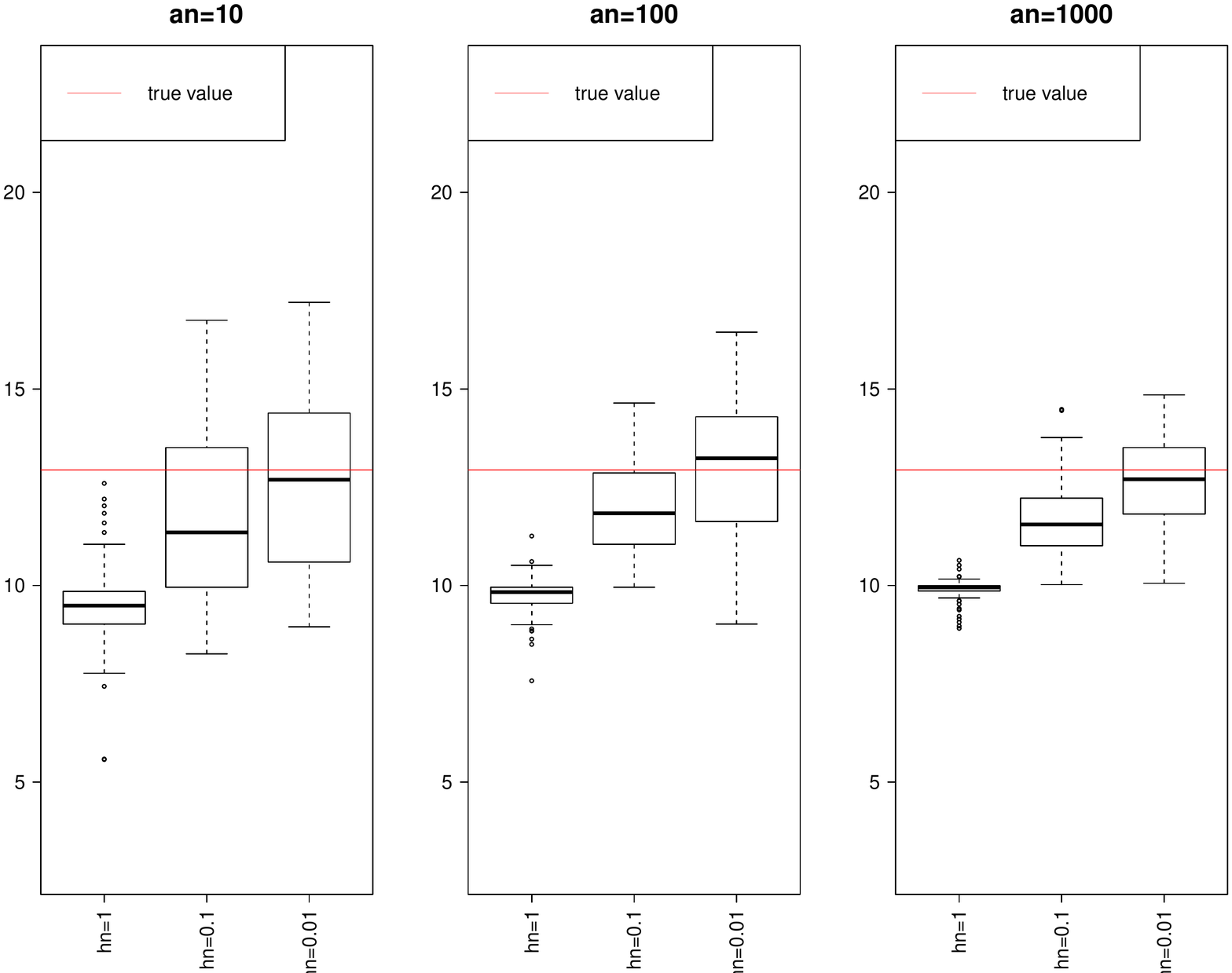}
  \end{center}
  \caption{Box plots for the estimated values of the optimal dividend barrier for $\alpha_n=10$ (left), $\alpha_n=100$ (middle), and $\alpha_n=1000$ (right). In each figure, $h_n=1$ (left), $h_n=0.1$ (middle), and $h_n=0.01$ (right). The red line shows the true value.} \label{fig5}
\end{figure}

Now, we examine mean, standard deviation (std), bias, and MSE for $h_n=1,0.1,0.01$ and $\alpha_n=10,100,1000$. 
We generate $100$ replications for each run of the simulations. 
Figure \ref{fig5} shows the box-plot for estimated values $\hat{\vartheta}_n^{(j)}$ for each $\alpha_n$ and $h_n$. 
The left figure is the case of $\alpha_n=10$, the middle figure is the case of $\alpha_n=100$, and the right figure is the case of $\alpha_n=1000$. 
In each figure, the left box is the case of $h_n=1$, the middle box is the case of $h_n=0.1$, and the right box is the case of $h_n=0.01$. 
The red line shows the true value $\vartheta_0=12.93958$. 
In view of the median (and mean) of the estimated values, it can be seen that the value converges to the true value as $\alpha_n$ increases and $h_n$ decreases. 
On the other hand, in view of the dispersion, it looks that the estimated values shrink as $\alpha_n$ increases. 
However, when $h_n$ is not sufficiently small, it seems that the estimated values converges to a value different from the true value. 
This phenomenon indicates that it is a warning that an asymptotic bias will occur unless the sampling interval $h_n$ is sufficiently small. 
Table 1 
shows mean $\left(\mu_n:=\frac{1}{B}\sum_{j=1}^B\hat{\vartheta}_n^{(j)}\right)$, std $\left(\sigma_n:=\sqrt{\frac{1}{B}\sum_{j=1}^B(\hat{\vartheta}_n^{(j)}-\vartheta_0)^2}\right)$, bias $\left(\frac{1}{B}\sum_{j=1}^B(\hat{\vartheta}_n^{(j)}-\vartheta_0)\right)$, and MSE $\left(\mu_n^2+\sigma_n^2\right)$ of the estimated values $\hat{\vartheta}_n^{(j)}$ for $\alpha_n=10,100,1000$, $h_n=1,0.1,0.01$ and $B=100$. 
It can be seen that the mean converges to the true value, and the MSE converges to $0$ as $\alpha_n$ increases and $h_n$ decreases. 
Note that the bias tends to have negative value which implies that the distribution of the estimator tends to be asymmetric. 
From the insurer's point of view, this phenomenon is a warning because setting a lower dividend barrier poses a risk to insures.      

\begin{table}[htbp]
\label{tab1}
 \caption{Mean, std, bias and MSE for the estimated values in 100 replications where the true optimal dividend barrier is $\vartheta_0 = 12.93958$.}
 \begin{center}
  \begin{tabular}{|c|c||c|c|c|c|}
    \hline
    $\alpha_n$   & $h_n$   &  mean  &  std  & bias   &  MSE  \\
    \hline
    \hline
              & 1      & 9.482    & 1.039   & -3.45   & 13.03\\
      10      & 0.1    & 11.778   & 2.106   & -1.16   & 5.78   \\
              & 0.01   & 12.560   & 2.126   & -0.37   & 4.66    \\
    \hline
              &  1     & 9.681    & 0.489   & -3.25   & 10.85  \\
      100     &  0.1   & 12.009   & 1.290   & -0.92   &  2.52  \\
              &  0.01  & 12.964   & 1.558   & 0.02   &  2.42  \\
    \hline
              & 1      & 9.895    & 0.260   & -3.04   & 9.33  \\
      1000    & 0.1    & 11.680   & 0.998   & -1.25   & 2.58   \\
              & 0.01   & 12.680   & 1.183   & -0.25   & 1.46   \\
    \hline
  \end{tabular}
 \end{center}
\end{table}

\section{Proofs}\label{proofs}
\noindent This section shows proofs of lemmas and theorems.
\subsection{Proof of Lemma 1}
From the definition, we can write
\begin{align*}
h^{\vartheta}(\hat{X}^{i,n})= \int_0^{\tau^{\vartheta}(\hat{X}^{i,n})}e^{-rt}d\xi_t^{\vartheta}(\hat{X}^{i,n})
&\le\int_0^{\infty} e^{-rt}d \overline{\hat{X}}^{i,n}_t \\
&= \sum_{k=1}^n e^{-rt_k}(\Delta_{i(k)}X \vee 0)\le \sum_{k=1}^n e^{-rt_k}|\Delta_{i(k)}X|,
\end{align*}
where $\overline{\hat{X}}^{i,n}_t= \sup_{0\le s<t}\hat{X}^{i,n}_s$ and $\{|\Delta_{i(k)}X|\}_{k=1,\ldots,n}$ is a sequence of i.i.d. random variables. 
From (\ref{eq1}), we have
$$ |\Delta_{i(k)}X| \stackrel{d}{=}|\Delta_{k}X| = |X_k-X_{k-1}| \stackrel{d}{=} |ch_n +\sigma W_{h_n}-S_{h_n}|\le ch_n +\sigma |W_{h_n}|+ S_{h_n}.$$ 
It is easy to see $\mathbb{E}\left[|W_{h_n}|^k\right]\lesssim h_n$. 
Since $\varphi_n(\lambda):= \mathbb{E}\left[e^{i\lambda S_{h_n}}\right] = e^{h_n \psi_S(\lambda)}$, we have
\begin{align*}
\mathbb{E}\left[S_{h_n}\right] &= i^{-1}\frac{d \varphi_n(\lambda)}{d\lambda}\bigr|_{\lambda=0}=h_n\mathbb{E}[S_1] \lesssim h_n, \\
\mathbb{E}\left[S_{h_n}^2\right] &= - \frac{d^2 \varphi_n(\lambda)}{(d\lambda)^2}\bigr|_{\lambda=0} =h_n\mathbb{E}[S_1^2]-h_n^2\mathbb{E}[S_1]^2 \lesssim h_n, 
\end{align*}
which imply that $\mathbb{E}\left[|\Delta_{i(k)}X|^m\right]\lesssim h_n$ for $m=1,2$. 
By the Taylor expansion, we have 
$$ \sum_{k=1}^n e^{-rt_k}= \frac{e^{-rh_n}(1-e^{-rnh_n})}{1-e^{-rh_n}} \lesssim \frac{1}{1-e^{-rh_n}} \lesssim h_n^{-1}.$$
Therefore, we have 
\begin{align*}
\mathbb{E}\left[h^{\vartheta}(\hat{X}^{i,n})\right] &\le \sum_{k=1}^n e^{-rt_k}\mathbb{E}\left[|\Delta_{i(k)}X|\right] =O(h_n^{-1}) O(h_n)=O(1),\\
\mathbb{E}\left[h^{\vartheta}(\hat{X}^{i,n})^2 \right] &\le \sum_{k=1}^n e^{-2rt_k}\mathbb{E}\left[|\Delta_{i(k)}X|^2\right] =O(h_n^{-1}) O(h_n)=O(1). \qed
\end{align*}

\subsection{Proof of Lemma 2}   
By definition, $\tau_n^{\vartheta}\in \{t_1,t_2,\ldots,t_n\}$ for any $X\in\tilde{D}_{\infty}$ and $\vartheta\in\Theta$. 
Note that if $U_{t_n}^{\vartheta}(\hat{X}^{i_{\mathrm{id}},n}):=\hat{X}_{t_n}^{i_{\mathrm{id}},n}-\xi_{t_n}^{\vartheta}(\hat{X}^{i_{\mathrm{id}},n})\ge 0$, we define $\tau^{\vartheta}(\hat{X}^{i_{\mathrm{id}},n})=\tau_n^{\vartheta}(X)=t_n$. 
This implies that we can divide $\mathcal{H}_n$ into $\mathcal{H}_{n,k}$ ($k=1,\ldots,n$), where $\mathcal{H}_{n,k}=\{h_{n,k}^{\vartheta},\vartheta\in\Theta\}$ with 
\begin{align*}
h_{n,k}^{\vartheta}(X)= \int^{t_k}_{0}e^{-rt}d\xi_{n,t}^{\vartheta}(X).
\end{align*}
For a fixed $X\in\tilde{\mathbb{D}}_{\infty}$, it can be seen that $\xi_{n,t}^{\vartheta}(X)\le \xi_{n,t}^{\vartheta'}(X)$ for any $n\in\mathbb{N},t>0$ if $\vartheta\ge \vartheta'$, which implies that
$$ \mathbb{P}_{\alpha_n}\left(\bm{1}_{\left\{h_{n,k}^{\bar{\vartheta}}(X) \le h_{n,k}^{\vartheta}(X)\le h_{n,k}^{u}(X)\right\}}\right) =1,\quad \forall\vartheta\in \Theta=(u,\bar{\vartheta}).$$
In addition, since $u\le \overline{\hat{X}}_{t_{l-1}}^{i_{\mathrm{id}},n}
\le \overline{\hat{X}}_{t_l}^{i_{\mathrm{id}},n}
$ ($l=1,\ldots,n$), where $\overline{\hat{X}}_{t}^{i_{\mathrm{id}},n}= \sup_{0\le s<t}\hat{X}^{i_{\mathrm{id}},n}_s$, we can write
\begin{align*}
&\left|\left\{\xi_{n,t_{l}}^{u}(\hat{X}^{i_{\mathrm{id}},n})-\xi_{n,t_{l-1}}^{u}(\hat{X}^{i_{\mathrm{id}},n})\right\}-\left\{\xi_{n,t_{l}}^{\bar{\vartheta}}(\hat{X}^{i_{\mathrm{id}},n})-\xi_{n,t_{l-1}}^{\bar{\vartheta}}(\hat{X}^{i_{\mathrm{id}},n})\right\}\right|\\ 
&= \left|\left(u\vee \overline{\hat{X}}_{t_{l-1}}^{i_{\mathrm{id}},n}\right)-\left(u\vee \overline{\hat{X}}_{t_{l}}^{i_{\mathrm{id}},n}\right)-\left(\bar{\vartheta}\vee \overline{\hat{X}}_{t_{l-1}}^{i_{\mathrm{id}},n}\right)+\left(\bar{\vartheta}\vee \overline{\hat{X}}_{t_{l}}^{i_{\mathrm{id}},n}\right)\right|\\
&= \bm{1}_{\{\overline{\hat{X}}_{t_{l-1}}^{i_{\mathrm{id}},n} \le \bar{\vartheta}\}}\left|\left(\overline{\hat{X}}_{t_{l}}^{i_{\mathrm{id}},n} -\overline{\hat{X}}_{t_{l-1}}^{i_{\mathrm{id}},n} \right) \wedge \left(\bar{\vartheta}- \overline{\hat{X}}_{t_{l-1}}^{i_{\mathrm{id}},n}\right)\right|\\
&\le \left|\Delta_l X \right|,
\end{align*}
which implies that 
\begin{align*}
\|h_{n,k}^{u}-h_{n,k}^{\bar{\vartheta}}\|_{\mathbb{P}_{\alpha_n},1}
&=\mathbb{P}_{\alpha_n}\left(|h_{n,k}^{u}(X)-h_{n,k}^{\bar{\vartheta}}(X)|\right)\\
&=\mathbb{P}_{\alpha_n}\left(\left| \int^{t_k}_{0} e^{-rt}\left\{d\xi_{n,t}^{u}(X)-d\xi_{n,t}^{\bar{\vartheta}}(X)\right\}\right|\right)\\
&\le \sum_{l=1}^k e^{-rt_l}\mathbb{P}_{\alpha_n}\biggl(  \biggl|\left\{\xi_{n,t_{l}}^{u}(\hat{X}^{i_{\mathrm{id}},n})-\xi_{n,t_{l}}^{u}(\hat{X}^{i_{\mathrm{id}},n})\right\}\\
& \hspace{20mm}   -\left\{\xi_{n,t_{l}}^{\bar{\vartheta}}(\hat{X}^{i_{\mathrm{id}},n})-\xi_{n,t_{l}}^{\bar{\vartheta}}(\hat{X}^{i_{\mathrm{id}},n})\right\}\biggr| \biggr)\\
&\le \sum_{l=1}^k e^{-rt_l}\mathbb{P}_{\alpha_n}\left(\left|\Delta_l X\right|\right)\\
&\stackrel{d}{=} \left(\sum_{l=1}^k e^{-rt_l}\right)\left|\Delta_1 X\right|.
\end{align*}
Since $\sum_{l=1}^k e^{-rt_l}=O(h_n^{-1})$ and $\mathbb{E}\left[\left|\Delta_1 X\right|\right]=O(h_n)$, we have for any $\epsilon>0$
\begin{align*}
\sum_{k=1}^n \mathbb{E}\left[\frac{\|h_{n,k}^{u}-h_{n,k}^{\bar{\vartheta}}\|_{\mathbb{P}_{\alpha_n},1}}{\epsilon}\right] 
&\le \frac{1}{\epsilon }\sum_{k=1}^n \left(\sum_{l=1}^k e^{-rt_l}\right) \mathbb{E}\left[\left|\Delta_{1}X\right|\right] \\
&=\frac{n}{\epsilon }O(h_n^{-1})O(h_n)=O(n).
\end{align*}
Note that the $L^1(\mathbb{P}_{\alpha_n})$-size of the brackets is bounded by $\epsilon$, which imples that  
\begin{align*}
\mathbb{E}\left[N_{[]}(\epsilon,\mathcal{H}_{n},L^1(\mathbb{P}_{\alpha_n}))\right]
&=\sum_{k=1}^n \mathbb{E}\left[N_{[]}(\epsilon,\mathcal{H}_{n,k},L^1(\mathbb{P}_{\alpha_n}))\right] \\
&\le \sum_{k=1}^n \mathbb{E}\left[\frac{\|h_{n,k}^{u}-h_{n,k}^{\bar{\vartheta}}\|_{\mathbb{P}_{\alpha_n},1}}{\epsilon}+1\right].
\end{align*}
Therefore, from Assumption 2, it follows that
\begin{align*}
\mathbb{E}\left[N_{[]}(\epsilon,\mathcal{H}_{n},L^1(\mathbb{P}_{\alpha_n}))\right]=o(\alpha_n).
\end{align*}
From the relationship between bracketing number and covering number (cf., Kosorok (2008; Lemma 9.18)
), we have 
\begin{align*}
\mathbb{E}\left[N(\epsilon,\mathcal{H}_{n},L^1(\mathbb{P}_{\alpha_n}))\right]\le \mathbb{E}\left[N_{[]}(\epsilon,\mathcal{H}_{n},L^1(\mathbb{P}_{\alpha_n}))\right]=o(\alpha_n).\qed
\end{align*}

\subsection{Proof of Lemma 3}
By the symmetrization result (cf., Kosorok (2008; Theorem 8.8)
), we can write 
\begin{align*}
\mathbb{E} \left[\|\mathbb{P}_{\alpha_n}-P\|_{\mathcal{H}_n}\right] 
&\le 2\mathbb{E}_X\left[\mathbb{E}_{\epsilon}\left[ \sup_{\vartheta\in \Theta}\left|\mathbb{P}_{\alpha_n} \left(\epsilon h_n^{\vartheta}(X)\right)\right|\Bigr| X \right]  \right]\\
&= 2\mathbb{E}_X\left[\mathbb{E}_{\epsilon}\left[ \sup_{\vartheta\in \Theta}\left|\frac{1}{\alpha_n}\sum_{i\in A_n} \epsilon^{(i)} h_n^{\vartheta}(X^{(i)})\right| \Bigr| X^{(i)}, i\in A_n \right]\right],
\end{align*}
where $\{\epsilon^{(i)}\}$ is a sequence of independent Rademacher random variables which are independent of $\{X^{(i)}\}$ and satisfy $\mathbb{P}(\epsilon^{(i)}=-1)=\mathbb{P}(\epsilon^{(i)}=1)=1/2$, and $\mathbb{E}_X,\mathbb{E}_{\epsilon}$ are the expectations with respect to $X^{(i)},\epsilon^{(i)}$, respectively. 
\color{black}
For any fixed $n\in\mathbb{N}$, $\delta>0$ and $\{X^{(i)}\}_{i\in A_n}$, let $\mathcal{H}_{n,j}(j=1,\ldots,N(\delta,\mathcal{H}_n,L^1(\mathbb{P}_{\alpha_n})))$ be 
a sequence of finite $\delta$-balls in $L^1(\mathbb P_{\alpha_n})$ over $\mathcal H_n$ (i.e., $\mathcal{H}_{n,j}$ is a subset of $\mathcal{H}_n$ and for any $h_n^{\vartheta},h_n^{\vartheta'}\in \mathcal{H}_{n,j}$, $\|h_n^{\vartheta}-h_n^{\vartheta'}\|_{\mathbb P_{\alpha_n},1}< \delta$ and $\cup_{j}\mathcal{H}_{n,j} \supset \mathcal{H}_n$).  
For each $\mathcal{H}_{n,j}$, we fix $\vartheta_j$ (satisfying   $\vartheta_j\neq \vartheta_{j'}$ if $j\neq j'$) which is a representative $h_n^{\vartheta_j}$ such that for any $h_n^{\vartheta}\in\mathcal{H}_{n,j}$
\begin{align*}
&\mathbb{E}_{\epsilon}\left[ \left|\frac{1}{\alpha_n}\sum_{i\in A_n} \epsilon^{(i)} h_n^{\vartheta}(X^{(i)})\right| \Bigr| X^{(i)}, i\in A_n \right]\\
&\le \mathbb{E}_{\epsilon}\left[ \left|\frac{1}{\alpha_n}\sum_{i\in A_n} \epsilon^{(i)} h_n^{\vartheta_j}(X^{(i)})\right| \Bigr| X^{(i)}, i\in A_n \right] +
\delta, 
\end{align*}
\color{black}
which implies that
\begin{align}
&\mathbb{E}_{\epsilon} \left[\sup_{\vartheta\in \Theta}\left|\frac{1}{\alpha_n}\sum_{i\in A_n} \epsilon^{(i)} h_n^{\vartheta}(X^{(i)})\right| \Bigr| X^{(i)}, i\in A_n \right]\nonumber\\
&\le 
 \mathbb{E}_{\epsilon}\left[\max_{j} \left|\frac{1}{\alpha_n}\sum_{i\in A_n} \epsilon^{(i)} h_n^{\vartheta_j}(X^{(i)})\right| \Bigr| X^{(i)}, i\in A_n \right]+\delta.\label{eq9}
\end{align}
\color{black}
Let $Z_j=\frac{1}{\alpha_n}\sum_{i\in A_n} \epsilon^{(i)} h_n^{\vartheta_j}(X^{(i)})$ and $Z=\max_j|Z_j|$. 
Then,
$$\|Z\|_{1|X}:=\mathbb{E}_{\epsilon}[|Z|\ | X^{(i)}, i\in A_n]\le \mathbb{E}_{\epsilon}[Z^2| X^{(i)}, i\in A_n]^{1/2}=:\|Z\|_{2|X},$$ 
from Jensen's inequality. 
On the other hand, based on the nondecreasing, nonzero convex function $\psi_2(x)=\exp(x^2)-1$, we introduce the Orlicz-norm 
$$\|Z\|_{\psi_2|X}:=\inf\left\{c>0 \Bigr|\ \mathbb{E}_{\epsilon}\left[\frac{\psi_2(|Z|)}{c} \Bigr|X^{(i)}, i\in A_n\right] \le 1 \right\},$$
for which $\|Z\|_{2|X} \le \|Z\|_{\psi_2|X}$. 
Applying the maximal inequality (cf., Kosorok (2008; Lemma 8.2)
), we have 
\begin{align*}
\|Z\|_{\psi_2|X} = \|\max_j|Z_j|\|_{\psi_2|X} \le K \psi_2^{-1}(N(\delta,\mathcal{H}_n,L^1(\mathbb{P}_{\alpha_n}))) \max_{j} \|Z_j\|_{\psi_2|X},
\end{align*}
where the constant $K$ depends only on $\psi_2$, which implies that the left-hand-side of (\ref{eq9}) is bounded by 
\begin{align*}
\sqrt{\log \left\{1+ N(\delta,\mathcal{H}_n,L^1(\mathbb{P}_{\alpha_n}))\right\}} \max_{j} \left\|Z_j\right\|_{\psi_2|X}+\delta,
\end{align*}
\color{black}
up to a constant. By Hoeffiding's inequality (cf., Kosorok (2008; Lemma 8.7)
), we have 
\begin{align*}
\mathbb{E}_{\epsilon}\left[\bm{1}_{\{|Z_j|>x\}} \Bigr|X^{(i)}, i\in A_n\right]
=\mathbb{P}_{\epsilon}\left(|Z_j|>x \Bigr|X^{(i)}, i\in A_n \right)
\le 2\exp\left(-\frac{1}{2}x^2/\|Z_j\|_{\epsilon}^{2}\right),
\end{align*}
for any $x>0$ and each $j$, where $\|Z_j\|_{\epsilon} = \mathbb{E}_{\epsilon}[|Z_j| \ |X^{(i)}, i\in A_n]$. 
Hence, from Kosorok (2008; Lemma 8.1) 
and Jensen's inequality, 
\begin{align*}
\left\|Z_j\right\|_{\psi_2|X}
&\le \left(\frac{1+2}{1/(2\|Z_j\|_{\epsilon}^{2})}\right)^{1/2}\\
&= \sqrt{6}\|Z_j\|_{\epsilon}\\
&\le \frac{\sqrt{6}}{\alpha_n}\mathbb{E}_{\epsilon}\left[\left|\sum_{i\in A_n}\epsilon^{(i)}h_n^{\vartheta_j}(X^{(i)})\right| \Bigr|X^{(i)}, i\in A_n\right]   \\
&\le \frac{\sqrt{6}}{\alpha_n} \left\{ \mathbb{E}_{\epsilon}\left[\left|\sum_{i\in A_n}\epsilon^{(i)}h_n^{\vartheta_j}(X^{(i)})\right|^2 \Bigr|X^{(i)}, i\in A_n\right]\right\}^{1/2}\\
&= \sqrt{\frac{6}{\alpha_n}} \left\{ \frac{1}{\alpha_n}\sum_{i\in A_n}\left| h_n^{\vartheta_j}(X^{(i)})\right|^2\right\}^{1/2}\\
&= \sqrt{\frac{6}{\alpha_n}}\sqrt{\mathbb{P}_{\alpha_n} \left((h_n^{\vartheta_j})^2\right)}, 
\end{align*}
which, together with Lemma 1, implies that
\begin{align*}
\mathbb{E}_X\left[\left\|Z_j\right\|_{\psi_2|X}\right] 
&\le \sqrt{\frac{6}{\alpha_n}}  \mathbb{E}_X \left[\sqrt{\mathbb{P}_{\alpha_n}\left((h_n^{\vartheta_j})^2\right)}\right] \\
&\le \sqrt{\frac{6}{\alpha_n}}  \left\{\mathbb{E}_X \left[\mathbb{P}_{\alpha_n}\left((h_n^{\vartheta_j})^2\right)\right]\right\}^{1/2}  
= O\left(\frac{1}{\sqrt{\alpha_n}}\right),
\end{align*}
uniformly in $\vartheta_j\in\Theta$. 
From this and Lemma 2, $\mathbb{E} \|\mathbb{P}_{\alpha_n}-P\|_{\mathcal{H}_n} $ converges to $0$ as $n\to\infty$ and $\delta\to 0$. \qed

\subsection{Proof of Lemma 4}
From the definition, we can write for any $\vartheta \in\Theta$
\begin{align*}
P(h_n^{\vartheta}-h^{\vartheta}) 
&= \mathbb{E}\left[\int^{\tau_n^{\vartheta}(X)}_{0} e^{-rt}d\xi_{n,t}^{\vartheta}(X) \right] - \mathbb{E}\left[\int^{\tau^{\vartheta}(X)}_{0} e^{-rt}d\xi_{t}^{\vartheta}(X) \right]\\
&= \mathbb{E}\left[\int^{\infty}_{0}\left\{\bm{1}_{\{\tau_n^{\vartheta}(X)> t\}}-\bm{1}_{\{\tau^{\vartheta}(X)> t\}}\right\}  e^{-rt}d\xi_{n,t}^{\vartheta}(X) \right] \\
&+ \mathbb{E}\left[\int^{\infty}_{0} \bm{1}_{\{\tau^{\vartheta}(X)> t\}}\left\{\bm{1}_{\{\overline{\hat{X}}_{t}^{i,n}> \vartheta\}}-\bm{1}_{\{\bar{X}_{t}> \vartheta\}}\right\}e^{-rt}d\overline{\hat{X}}_{t}^{i,n}\right]\\
&+ \mathbb{E}\left[\int^{\infty}_{0} \bm{1}_{\{\tau^{\vartheta}(X)> t\}}\bm{1}_{\{\bar{X}_{t} > \vartheta\}}e^{-rt}d\left\{ \overline{\hat{X}}_t^{i,n} - \bar{X}_t\right\} \right]\\
&=:I_1+I_2+I_3\quad (\mathrm{say}),
\end{align*}
where $\overline{\hat{X}}^{i,n}_t= \sup_{0\le s< t}\hat{X}_{s}^{i,n}$. 
For the term $I_1$, Lemma 1 yields
\begin{align}
|I_1| &\le \left|\mathbb{E}\left[\int^{T}_{0}\left\{\bm{1}_{\{\tau_n^{\vartheta}(X)> t\}}-\bm{1}_{\{\tau^{\vartheta}(X)> t\}}\right\} e^{-rt}d\xi_{n,t}^{\vartheta}\right]\right|+O\left(e^{-rT}\right), \label{eq10}
\end{align}
for any fixed $T>0$. 
Denoting $F_{n,\tau}^{\vartheta}(t)=\mathbb{P}\left(\tau_n^{\vartheta}(X)\le t \right)$ and $F_{\tau}^{\vartheta}(t)=\mathbb{P}\left(\tau^{\vartheta}(X)\le t \right)$, the first term of the right hand side of (\ref{eq10}) is bounded by 
\begin{align*}
\mathbb{E}\left[\int^{\infty}_{0} e^{-rt}d\xi_{n,t}^{\vartheta} \right]\left| F_{n,\tau}^{\vartheta}(T) - F_{\tau}^{\vartheta}(T)\right|.
\end{align*}
\color{black}
Let $\tilde{\mathbb{D}}_{T}\subset \tilde{\mathbb{D}}_{\infty}$ be a space of c\`{a}dl\`{a}g functions on $[0,T]$. 
We now consider the Skorokhod topology $(\tilde{\mathbb{D}}_{T},d_T)$, where $d_T$ is the Skorokhod metric defined by
$$ d_T(x,y)=\inf_{\lambda\in\Lambda_T}\left(\max\left\{\|x\circ \lambda-y\|_T,\|\lambda-I\|_T\right\}\right),$$
for any $x=(x_t),y=(y_t)\in \tilde{\mathbb{D}}_{T}$ (cf., \cite{B99}). 
Note that $\Lambda_T$ is the class of strictly increasing, continuous mappings of $[0,T]$ onto itself, $x\circ \lambda=(x_{\lambda_t})$ for any $\lambda=(\lambda_t)\in\Lambda_T$, $I$ is the identity map on $[0,T]$ and $\|z\|_T=\sup_{0<t\le T}z_t$. 
On this topology, we define a map $g^{\vartheta}:(\tilde{\mathbb{D}}_{T},d_T)\to (\bar{\mathbb{R}},|\cdot|)$ by
\begin{align*}
g^{\vartheta}(x) = \inf_{0<t\le T}\left\{ x_t - \left(\sup_{0<s<t}x_s-\vartheta \right)\vee 0 \right\}.
\end{align*}
Then, it is easy to see that
\begin{align*}
&|g^{\vartheta}(x)-g^{\vartheta}(y)| \\
&=\left|\inf_{0<t\le T}\left\{x_t + \left(\inf_{0<s<t}(-x_s)+\vartheta\right)\wedge 0\right\} - \inf_{0<t\le T}\left\{y_t + \left(\inf_{0<s<t}(-y_s)+\vartheta\right)\wedge 0\right\} \right|\\
&\le \sup_{0<t\le T}|\inf_{0<s<t}x_s-\inf_{0<s<t}y_s| \\
&\hspace{10mm} + \sup_{0<t\le T}\left|\left\{ \inf_{0<s<t}(-x_s) +\vartheta\right\}\wedge 0 - \left\{\inf_{0<s<t}(-y_s)+\vartheta\right\}\wedge 0\right|\\
&\le \sup_{0<t\le T}|x_t-y_t| + \sup_{0<t\le T}\left|\sup_{0<s<t}x_s - \sup_{0<s<t}y_s\right|\\
&\le 2\sup_{0<t\le T}|x_t-y_t|\\
&\lesssim d_T(x,y),
\end{align*}
which implies that $g$ is continuous on $(\tilde{\mathbb{D}}_{T},d_T)$. Therefore, by using the continuous mapping theorem, Theorem 1 implies
\begin{align*}
\inf_{0<t\le T}U_t^{\vartheta}(\hat{X}^{i_{\mathrm{id}},n}) = g^{\vartheta}(\hat{X}^{i_{\mathrm{id}},n}) \leadsto g^{\vartheta}(X)=\inf_{0<t\le T}U_t^{\vartheta}(X).   
\end{align*}
From the definition of the weak convergence, it follows that 
\begin{align*}
|F_{n,\tau}^{\vartheta}(T)-F_{\tau}^{\vartheta}(T) |
&= \left|\mathbb{P}\left(\tau_{n}^{\vartheta}\le T\right) - \mathbb{P}\left(\tau^{\vartheta}\le T\right)\right|\\
&= \left|\mathbb{P}\left(\inf_{0<t\le T}U_t^{\vartheta}(\hat{X}^{i_{\mathrm{id}},n})<0\right) - \mathbb{P}\left(\inf_{0<t\le T}U_t^{\vartheta}(X)<0 \right)\right|\to 0,  
\end{align*}
as $n\to\infty$ for any $\vartheta\in \Theta$. 
Since $\mathbb{E}\left[\int_0^{\infty}e^{-rt}d\xi_{n,t}^{\vartheta} \right]=O(1)$ from Lemma 1, we have $|I_1|\to 0$ as $n\to\infty$ and $T\to\infty$. 
In the same way, we have $|I_2|\to 0$ as $n\to\infty$ and $T\to\infty$.
\color{black}
For the term $I_3$, Lemma 1 yields
\begin{align*}
|I_3| \le \left| \mathbb{E}\left[\int^{T}_{0} \bm{1}_{\{\tau^{\vartheta}(X)> t\}}\bm{1}_{\{\bar{X}_{t} > \vartheta\}}e^{-rt}d\left\{ \overline{\hat{X}}_t^{i,n} - \bar{X}_t\right\}\right]\right|+O(e^{-rT}),
\end{align*}
for any fixed $T>0$. 
Then, there exits a constant $M>0$ such that 
\begin{align}
&\left| \mathbb{E}\left[\int^{T}_{0} \bm{1}_{\{\tau^{\vartheta}(X)> t\}}\bm{1}_{\{\bar{X}_{t} > \vartheta\}}e^{-rt}d\left\{ \overline{\hat{X}}_t^{i,n} - \bar{X}_t\right\}\right]\right|\nonumber \\
&\le \mathbb E \left[\sup_{t\in [0,T]} \frac{\left|\overline{\hat{X}}_t^{i,n} - \bar{X}_t\right|}{\bar{X}_t}  \int^{T}_{0} e^{-rt}d\bar{X}_t   \right]\nonumber\\
&\le M \mathbb E \left[\sup_{t\in [0,T]} \left|\hat{X}_t^{i,n} - X_t\right|  \int^{T}_{0} e^{-rt}d\bar{X}_t   \right] \label{eq12}\\
&\le M \mathbb E \left[\left(\sup_{t\in [0,T]} \left|\hat{X}_t^{i,n} - X_t\right|\right)^2\right]^{1/2}  \mathbb E \left[\left(\int^{T}_{0} e^{-rt}d\bar{X}_t \right)^2  \right]^{1/2} \to 0, \label{eq13}
\end{align}
as $n\to\infty$. 
Note that (\ref{eq12}) is shown by $\sup_{t}|\sup_{s<t}x_s -\sup_{s<t}y_s|\le \sup_{t}|x_s-y_s|$ and $\bar{X}_t\ge u$, and (\ref{eq13}) is shown by $\mathbb E \left[\left(\sup_{t\in [0,T]} \left|\hat{X}_t^{i,n} - X_t\right|\right)^2\right]=o(1)$ from Theorem 1 and $\mathbb E \left[\left(\int^{T}_{0} e^{-rt}d\bar{X}_t \right)^2  \right]=O(1)$.   
Hence, $|I_3|\to 0$ as $n\to\infty$ and $T\to\infty$.  
Therefore, we have 
\begin{align*}
|P(h_n^{\vartheta}-h^{\vartheta})|\le |I_1|+|I_2|+|I_3|\to 0,
\end{align*}
as $n\to\infty$ and $T\to\infty$, uniformly in $\vartheta\in\Theta$.\qed

\subsection{Proof of Theorem 2}
Lemmas 3 and 4 imply that
\begin{align*}
\sup_{\vartheta\in\bar{\Theta}}|\mathbb{P}_{\alpha_n}h_n^{\vartheta}-Ph^{\vartheta}|
\le \sup_{\vartheta\in\bar{\Theta}}|(\mathbb{P}_{\alpha_n}-P)h_n^{\vartheta}|
+\sup_{\vartheta\in\bar{\Theta}}|P(h_n^{\vartheta}-h^{\vartheta})|
\stackrel{p}{\to}0,
\end{align*}
as $n\to\infty$. 
Combining this and (\ref{eq8}), we immediately have the conclusion from van der Vaart (1998; Theorem 5.7)
.\qed

\begin{acknowledgement}
The autors would like to thank the Editors for their constructive comments. 
The first author was partially supported by JSPS KAKENHI Grant Number JP21K03358 and JST CREST JPMJCR14D7, Japan. 
The second author was supported by JSPS KAKENHI Grant Number JP21K11793. 
\end{acknowledgement}

\end{document}